\documentclass[a4paper,12pt]{article}
\usepackage{amssymb,amsmath}

\newcommand{\Z}{\mathbb{Z}}

\newcommand{\F}{\mathbb{F}}

\newcommand{\beql}[1]{\begin{equation}\label{#1}}
\newcommand{\eeq}{\end{equation}}

\newcommand{\modd}[1]{\; ( \text{mod} \; #1)}
\newtheorem{theorem}{Theorem}
\newtheorem{corollary}{Corollary}
\newtheorem{lemma}{Lemma}
\newtheorem{definition}{Definition}
\newtheorem{prop}{Proposition}

\begin{document}

\title{Iteration of Quadratic Polynomials Over Finite Fields}

\author{D.R. Heath-Brown\\Mathematical Institute, Oxford}

\date{}

\maketitle



\section{Introduction}

Let $f(X)\in\F_q[X]$ and define the iterates $f^j(X)$ by setting
$f^0(X)=X$ and $f^{j+1}(X)=f(f^j(X))$.  Let $m\in\F_q$, and consider
the sequence of values $f^0(m),f^1(m),f^2(m),\ldots$. Since the
field $\F_q$ is finite, the
sequence eventually recurs, and one enters a closed cycle.  We are
interested in the questions:- How long is it before one enters the
cycle? How long is the cycle? In general we can construct a
directed graph $\Gamma_f=\Gamma_f(\F_q)$, 
whose vertices are the elements $m$ of $\F_q$, and
with edges $(m,f(m))$. The trajectory $f^0(m),f^1(m),f^2(m),\ldots$
then consists of a pre-cyclic ``tail'', followed by a cycle.

Linear polynomials are easily handled.
When $f(X)=X+b$ one has $f^j(X)=X+bj$, so that if $b=0$ the cycles are
singleton sets, and if $b\not=0$ then $\Gamma_f$ is a union of cycles
of length $p$, the characteristic of the field. For linear polynomials
$f(X)=aX+b$ with $a\not=0,1$ one finds that
\[f^j(X)=a^j\{X+b(a-1)^{-1}\}-b(a-1)^{-1}.\]
Thus $\Gamma_f$ 
consists of cycles of length ord$(a)$ together with a cycle
$\{-b(a-10^{-1}\}$ of length 1.

The situation is much more interesting for higher degree
polynomials, and forms the basis for Pollard's famous ``Rho
Algorithm'' for integer factorization \cite{pollard}.  If one wishes
to factor $N$ the algorithm calculates successive iterates $f^j(m)$ and
$f^{2j}(m)$ modulo $N$, until one reaches a value for which ${\rm
  g.c.d.}(f^j(m)-f^{2j}(m),N)>1$. If this highest common factor is
different from $N$ then one has obtained a non-trivial factor of
$N$. When $p$ is a prime divisor of $N$, the sequence of iterates
modulo $p$ will
have an initial segment of length $t$ say, (the ``tail'' of the letter
rho) followed by a cycle of length $c$ say.  Thus $p\mid
f^j(m)-f^{2j}(m)$ when $j$ is the smallest multiple of $c$ for which
$j>t$. In particular the first such $j$ is at most $t+c$.  If $p'$ is some
other prime divisor of $N$ there will be a corresponding value $j'$
for which $p\mid f^{j'}(m)-f^{2j'}(m)$. Unless the two values $j$ and
$j'$ are the same, the method will produce a nontrivial divisor
${\rm g.c.d.}(f^j(m)-f^{2j}(m),N)$ of $N$. The efficiency of the
algorithm depends on $t$ and $c$ being small.

A crude probabilistic argument predicts that, over the field $\F_q$,
the sequence $f^0(m),f^1(m),f^2(m),\ldots$ is likely to complete a
cycle after roughly $O(q^{1/2})$ steps. This is a version of the
``Birthday Paradox''.  Specifically, if one imagines the sequence
as taking values in $\F_q$ independently and uniformly at random, then
the chance of having a repetition within $N$ steps, say, is
\[1-\prod_{j=0}^{N-1}\left(1-\frac{j}{q}\right)\]
and when $N$ is of order $\sqrt{q}$ this is roughly $1-\exp\{N^2/2q\}$.
Thus there is a positive probability of a repetition as soon as
$N\gg\sqrt{q}$.

Unfortunately there are examples in which this heuristic clearly
fails. Thus if $f(X)=X^2$ one has $f^j(m)=m^{2^j}$, and if $m$ has odd
order $r$ one gets a pure cycle of length $l$, where $l$ is the order of 2
modulo $r$.  Thus if $q$ is a prime of the shape $2r+1$,
with $r$ a prime for which 2 is a primitive root, then the cycle
length will be $r-1=(q-3)/2$ whenever $m$ has order $r$ modulo $q$.
While it is not known that infinitely many such primes $q$ exist it is
certainly conjectured to be so. Thus we will expect to get cycles of
length $\gg q$ for a positive proportion of initial values $m$.

A second example is provided by the polynomial $f(X)=X^2-2$. If
$m=a+a^{-1}$ for some $a\in\F_q$, then $f^j(m)=a^{2^j}+a^{-2^j}$, and
we have a situation similar to that described above.  If $q=2r+1$
with $r$ a prime for which 2 is a primitive root, then again we will
have cycles of
length $\gg q$ for a positive proportion of initial values $m$.

Thirdly one can consider polynomials of the shape $f(X)=X^3+c$, in the
case in which $q$ is a prime with $q\equiv 2\modd{3}$. Here one sees
that $f$ induces a permutation of $\F_q$, since $X^3=a$
has a unique solution in $\F_q$, for every $a\in\F_q$.  If $m\in\F_q$
is given, the trajectory $f^0(m),f^1(m),\ldots$ is therefore
completely cyclic, and our question merely concerns the length of the
cycle. However the proportion of permutations in the symmetric group
$S_q$ for which $m$ belongs to a cycle of given length $k$, is exactly
$q^{-1}$. Thus one might expect all cycle lengths to occur equally
often, and that one should get a cycle of length  at least $q/2$, say,
with probability around $1/2$.  The numerical evidence seems to
support this. For a given prime $p\equiv 2\modd{3}$ and every
$c=1,\ldots p-1$ we compute the length, $l(c,p)$ say, of the
cycle which starts at $m=0$. We then see for how many values of $c$
the scaled cycle length $p^{-1}l(c,p)$ falls into each of the intervals
$((k-1)/10,k/10]$, for $k=1,\ldots 10$. If the permutations induced
by the various polynomials $X^3+c$ were genuinely random we would expect
roughly the same number of scaled cycle lengths in each such interval. The data
for the first two primes $p\equiv 2\modd{3}$ beyond $10^5$ are presented in
Table 1. The figures appear to support the random permutation model well.

\begin{table}
  \begin{centering}
\begin{tabular}{r|r|r}
  \mbox{Prime}& 100019 & 100043 \\ \hline
$(0,\tfrac{1}{10}]$& 10030 & 9936\\
$(\tfrac{1}{10},\tfrac{2}{10}]$& 9944 & 9730 \\
$(\tfrac{2}{10},\tfrac{3}{10}]$& 9992 & 9976 \\
$(\tfrac{3}{10},\tfrac{4}{10}]$& 10122 & 10232 \\
$(\tfrac{4}{10},\tfrac{5}{10}]$& 10212 & 10034 \\
$(\tfrac{5}{10},\tfrac{6}{10}]$& 9830 & 10000 \\
$(\tfrac{6}{10},\tfrac{7}{10}]$& 9902 & 10086 \\
$(\tfrac{7}{10},\tfrac{8}{10}]$& 9904 & 10012 \\
$(\tfrac{8}{10},\tfrac{9}{10}]$& 10070 & 9946 \\
$(\tfrac{9}{10},1]$ & 10012 & 10090\\ \hline
\end{tabular} 
\caption{Distribution of scaled cycle lengths}
\end{centering}
\end{table}

The main goal of the present paper is to describe a quite different
theory for the
iterates of quadratic polynomials in odd characteristic, 
in which it is clear why the
anomalous cases above must be excluded.  In contrast to the situation
with $f(X)=X^3+1$, when
$f(X)=aX^2+bX+c$, the equation $f(x)=s$ typically has either 2
solutions or none at all, the latter case holding for roughly half the
possible values of $s$ (those for which $b^2-4a(c-s)$ is a non-square).
When $f(x)=s$ has two solutions $x=t_1$ and $x=t_2$, the equations
$f(x)=t_1$ and $f(x)=t_2$ will again typically have either 2 solutions
or none.  In this way, considering solutions of $f^r(x)=m$, one sees
that $\Gamma_f$ is potentially much more complicated
than a series of cycles.  

Our main result demonstrates this distinction clearly.
\begin{theorem}\label{T1}
Let $\F_q$ be a finite field of characteristic $p\not=2$, and let
$f(X)=aX^2+c\in\F_q[X]$ with $a\not=0$. Suppose that $f^i(0)\not=f^j(0)$ for
$0\le i<j\le r$.  Then
\beql{stst}
 \# f^r(\F_q)=\mu_r q+O(2^{4^r} \sqrt{q})
\eeq
uniformly in $a$ and $c$, where the constant $\mu_r$ is defined recursively
by taking $\mu_0=1$ and
\beql{rr}
\mu_{r+1}=\mu_r-\tfrac12\mu_r^2.
\eeq

Moreover we have $\mu_r\sim 2/r$ as $r\rightarrow\infty$.
\end{theorem}

At this point we should mention some closely related work. Shao
\cite[Theorem 1.6]{shao} handles the case $f(X)=X^2+1$ by a method
which generalizes readily to other quadratics.  The condition in his
theorem is stronger than ours (that $f^i(0)\not=f^j(0)$ for
$0\le i<j\le r$) but an examination of the proof shows that he only
needs something like our condition.  His result does not include an
explicit dependence on $r$.  Juul, Kurlberg, Madhu and Tucker
\cite{JKMT} handle general rational functions rather than restricting
to quadratic polynomials. Their emphasis is on the reductions of a
given rational function $\phi(X)\in\mathbb{Q}(X)$ modulo different
primes, but they show under quite general conditions that the sum of
all cycle lengths is $o(p)$ as $p\rightarrow\infty$. (See Corollary 2
below.) 

Before discussing the implications of the theorem, let us examine the
condition that $f^i(0)\not=f^j(0)$ for $0\le i<j\le r$.  The critical
points of a polynomial $f(X)$ are the roots $\xi$ of $f'(X)$, and $f$
is said to be ``post-critically finite'' if the iterates $f^j(\xi)$
eventually enter a cycle, for every critical point $\xi$. In dynamics
in general post-critically finite maps are a very important
subclass.  Of course,
over a finite field every polynomial is post-critically finite.  However
our condition can be viewed as saying that, in an approximate sense,
$f$ fails to be post-critically finite. (When
$f(X)=aX^2+c$ the only critical point is $\xi=0$.)

Certainly the condition that $f^i(0)\not=f^j(0)$ for $0\le i<j\le r$
fails for the polynomials $f(X)=X^2$ and $f(X)=X^2-2$, with
$i=0,j=1$ and $i=2,j=3$ respectively. Suppose next that 
$f$ is the reduction of a
polynomial $F(X)=AX^2+C\in\Z[X]$, with $A, C>0$, then
the sequence $F^0(0),F^1(0),F^2(0),\ldots$ is strictly
monotonic, with $F^j(0)\le (A+C)^{2^j-1}$. Thus if $p\ge (A+C)^{2^r}$
we cannot have $p\mid F^j(0)-F^i(0)$ with $0\le i<j\le r$.  The
condition of the theorem will therefore hold when
\beql{rb}
r\le\frac{\log\log p}{\log 2}-\frac{\log\log(A+C)}{\log 2}.
\eeq
In following this paper the reader may wish to bear in mind the 
archetypal example
$f(X)=X^2+1$, for which $r\le\tfrac12\log\log p$ suffices.

Our main theorem above has the following immediate consequences.
\begin{corollary}
Let $\F_q$ be a finite field of characteristic $p\not=2$, and let
$f(X)=aX^2+c\in\F_q[X]$ with $a\not=0$. Then $f^i(0)=f^j(0)$ for some
$i,j$ with
\[i<j\ll \frac{q}{\log\log q}.\]
\end{corollary}

\begin{corollary}
Let $\F_p$ be a finite field with $p>2$ prime, and let
$f(X)=aX^2+c\in\F_q[X]$ be the reduction of $AX^2+C\in\Z[X]$, where
$A,C>0$.  Then the sum of all the cycle lengths in $\Gamma_f$ will be 
$O_{A,C}(p(\log\log p)^{-1})$.  Similarly the length of
any pre-cyclic path in $\Gamma_f$ will be $O_{A,C}(p(\log\log p)^{-1})$.
\end{corollary}

The first corollary gives an unconditional bound $o(q)$ for the first
recurrence in the sequence $f^0(0),f^1(0),f^2(0),\ldots$.  The second
corollary proves a similar result for arbitrary initial values for
the reductions of fixed positive definite quadratic polynomials
$AX^2+C$. Moreover it highlights the difference in behaviour between
such polynomials and the cubic case $f(X)=X^3+1$, where the cycle
lengths can sum to $q$.

To prove Corollary 1 we choose $r=[(\log\log q)/(\log 4)]-1$,
so that $2^{4^r}\sqrt{q}\ll q/r$. Then, according to Theorem 1, we
have either
$f^i(0)=f^j(0)$ for some $i<j\le r$, or $\# f^r(\F_q)\ll q/r$. Writing
the latter bound as $\# f^r(\F_q)\le Cq/r$ for an appropriate constant
$C$ we deduce in the latter case that if $k=[Cq/r]$ then the values
$f^r(0),f^{r+1}(0),\ldots,f^{r+k}(0)$ cannot be distinct, since they
all lie in $f^r(\F_q)$ and $k+1>Cq/r$. In either case there must
therefore be acceptable values $i<j\le r+k$.  The claim then follows.

For Corollary 2 we observe as above that the condition of the theorem
holds under the assumption (\ref{rb}). The choice
$r=[(\log\log  p)/(\log 4)]-1$, will
satisfy (\ref{rb}) when $p\gg_{A,C}1$, and the theorem then yields
$\# f^r(\F_p)\ll p/r$. All cycles lie inside $f^r(\F_p)$, giving the
first assertion of the corollary.  Moreover if $f^0(m),\ldots,f^t(m)$
is a pre-cyclic path then $f^r(m),\ldots,f^t(m)$ are distinct elements
in $f^r(\F_p)$, so that $t-r\ll p/r$. We then see that $t\ll r+p/r$,
from which the second assertion follows.

We should explain the restriction to polynomials $aX^2+c$. For an
arbitrary polynomial $f$, if we define $g(X):=f(X+d)-d$, then we will
have $g^j(X)=f^j(X+d)-d$.  Thus $\Gamma_g$ may be obtained from
$\Gamma_f$ by relabelling each vertex $m$ as $m-d$. Since the two
graphs are isomorphic in this sense, it suffices to study $f(X+d)-d$
for a suitably chosen $d$.  In the case in which $f(X)=aX^2+bX+c$ (and
$\F_q$ has odd characteristic) we can choose $d=-b/(2a)$ to produce a
polynomial $g(X)$ of the shape $aX^2+c'$. Thus we may translate our
results into statements about general quadratic polynomials as follows.

\begin{corollary}
Let $\F_q$ be a finite field of characteristic $p\not=2$, and let
$f(X)=aX^2+bX+c\in\F_q[X]$ with $a\not=0$. Suppose that 
$f^i(-b/(2a))\not=f^j(-b/(2a))$ for
$0\le i<j\le r$.  Then
\[ \# f^r(\F_q)=\mu_r q+O(2^{4^r} \sqrt{q})\]
uniformly in $a,b$ and $c$, with the same $\mu_r$ as before.

In particular $f^i(-b/(2a))=f^j(-b/(2a))$ for some
$i,j$ with
\[i<j\ll \frac{q}{\log\log q}.\]

If $q$ is prime, and $f$ is the reduction of a positive definite
quadratic polynomial 
$AX^2+BX+C\in\Z[X]$,  then the sum of all the cycle lengths in
$\Gamma_f$ will be  
$O_{A,B,C}(q(\log\log q)^{-1})$.  Similarly the length of
any pre-cyclic path in $\Gamma_f$ will be $O_{A,B,C}(q(\log\log q)^{-1})$.
\end{corollary}
In much the same way one can show that it would suffice to prove our
theorem for polynomials $f(X)=X^2+d$.  One could then deduce the
corresponding result for $aX^2+d/a$ by considering iterates of
$g(X):=a^{-1}f(aX)$.

Theorem 1 gives us an asymptotic formula $\# f^r(\F_q)\sim\mu_r
q$. We proceed to give a probabilistic argument showing why one might
expect this, and how the recurrence relation (\ref{rr}) arises.
When $r=0$ we have $\# f^0(\F_q)=q$, so that $\mu_0=1$. Suppose now that
we have a relation $\# f^r(\F_q)\sim\mu_r q$. We will use an inductive
argument to produce the corresponding result for $f^{r+1}$.

To have $m\in
f^{r+1}(\F_q)$ it is necessary and sufficient that $m\in f(\F_q)$ and
that $n\in f^r(\F_q)$ for at least one solution $n$ of $f(x)=m$.
Since $\F_q$ contains $(q+1)/2$ squares one has $m\in f(\F_q)$ in
exactly $(q+1)/2$ cases, and except for the value $m=f(0)$ there will
then be precisely two possible values of $n$.  Let these be $n_1$ and
$n_2$. If the probability of these lying in $f^r(\F_q)$ were $\mu_q$ each,
independently, one might expect that the probability of at least one
being in $f^r(\F_q)$ should be $2\mu_q-\mu_q^2$, by the
inclusion-exclusion principle. It would then follow that $m$ belongs
to $f^{r+1}(\F_q)$ with probability around $\tfrac12(2\mu_q-\mu_q^2)$.  One
would therefore produce an asymptotic expression $\# f^{r+1}(\F_q)\sim
\mu_{r+1}q$ with $\mu_{r+1}$ as in (\ref{rr}).

We next explain why $\mu_r\sim 2r^{-1}$, as claimed in Theorem
1. Writing $\nu_r=2/\mu_r$ we see that $\nu_0=2$ and
\[\nu_{r+1}=\nu_r+1+\frac{1}{\nu_r-1}.\]
An easy induction then shows that $\nu_r\ge r+2$ for all $r\ge 0$,
whence $\nu_{r+1}\le\nu_r+1+1/(r+1)$. Another induction shows that
\[\nu_r\le r+2+\sum_{j=1}^r j^{-1},\;\;\;(r\ge 1),\]
so that $\nu_r\le r+3+\log r$ for $r\ge 1$.  Together with the lower
bound $\nu_r\ge r+2$ this shows that $\nu_r\sim r$ and hence
$\mu_r\sim 2/r$.

{\bf Acknowledgments} The author would particularly like to extend his thanks
to Giacomo Micheli, for a number of interesting conversations introducing
the author to the subject of polynomial iteration. Joe Silverman also
provided a number of helful comments.  Thanks are also due
to Tim Browning, for elucidating a technical point in Section 3, to
Maksym Radziwi\l\l\, for some preliminary computational results, and to
Ben Green, Rafe Jones, Tom Tucker and Michael Zieve for some useful references.

\section{A Second Moment Calculation}

Fundamental to our treatment of Theorem 1 will be moments of the
functions
\[\rho_r(m)=\#\{x\in\F_q:f^r(x)=m\}.\]
Our first task is to estimate the moments
\[N(r;k):=\sum_{m\in \F_q}\rho_r(m)^k\]
for $r=0,1,2,\ldots$ and
$k=1,2,\ldots$. Trivially we have $\rho_0(m)=1$ for all $m$ so
that $N(0;k)=q$ for every $k$. Moreover it is also clear that
$N(r;1)=q$ for every $r$.

Before moving to the general situation it may be helpful to think
first about the case $k=2$, for which
\beql{N2}
N(r;2)=\#\{(x,y)\in\F_q^2:f^r(x)=f^r(y)\}.
\eeq
The equation $f^r(X)-f^r(Y)=0$ defines a curve in $\mathbb{A}^2$.  
An absolutely irreducible curve $C$ over $\F_q$ will have
$q+O_C(\sqrt{q})$ points, by Weil's ``Riemann Hypothesis''.  However our
curve is far from being irreducible.

Indeed
\[f^r(X)-f^r(Y)=
\left(f^{r-1}(X)+f^{r-1}(Y)\right)\left(f^{r-1}(X)-f^{r-1}(Y)\right),\]
whence a trivial induction produces
\beql{fact}
f^r(X)-f^r(Y)=(X-Y)\prod_{j=0}^{r-1}\left(f^j(X)+f^j(Y)\right).
\eeq
Thus we obtain $r+1$ factors.  However it is not immediately clear when
polynomials of the form $f^j(X)+f^j(Y)$ are absolutely irreducible
over $\F_q$.

In general, suppose that $\phi(X,Y)$ is a polynomial of degree $D$,
over a field
$K$, and let $\Phi(U,V,W)=W^D\phi(U/W,V/W)$ be the corresponding
form.  If $\Phi$ factors as $\Phi_1\Phi_2$ over the algebraic
completion $\overline{K}$ then there will necessarily be triple
$(u,v,w)\not=(0,0,0)\in\overline{K}^3$ such that
$\Phi_1(u,v,w)=\Phi_2(u,v,w)=0$.
For any such triple we then have $\nabla\Phi=\Phi_1\nabla\Phi_2+
\Phi_2\nabla\Phi_1=\mathbf{0}$. This gives us a simple criterion for
absolute irreducibility, which is sufficient, though not necessary: If
$\nabla\Phi$ vanishes only at the origin in $\overline{K}^3$, then
$\Phi$ must be absolutely irreducible.

We apply this criterion to $f^j(X)+f^j(Y)$. Writing $D=2^{j}$ for
convenience, and
\beql{nn}
F^j(U,W)=W^Df^j(U/W),
\eeq
we have
\begin{eqnarray*}
\lefteqn{\nabla(F^j(U,W)+F^j(V,W))}\\
&=&\left(W^{D-1}(f^j)'(U/W)\,,\,W^{D-1}(f^j)'(V/W)\,,\,
\frac{\partial}{\partial W}(F^j(U,W)+F^j(V,W))\right).
\end{eqnarray*}
If $f(X)=aX^2+c$ then $(f^j)'(X)=2af^{j-1}(X)(f^{j-1})'(X)$.  It then
follows by induction that
\[W^{D-1}(f^j)'(U/W)=(2a)^j\prod_{s=0}^{j-1}F^s(U,W).\]
In particular, if $\nabla(F^j(u,w)+F^j(v,w))$ vanishes, then there are
indices $s,t\le j-1$ for which $F^s(u,w)=F^t(v,w)=0$. Since
\[F^s(u,0)=a^{2^s-1}u^{2^s}\;\;\;\mbox{and}\;\;\; F^t(v,0)=a^{2^t-1}v^{2^t}\]
we see that $w=0$ would imply $u=v=w=0$, which is excluded.  We then
see that we would have $f^s(x)=f^t(y)=0$ for some $x,y\in
\overline{K}$ such that $f^j(x)+f^j(y)=0$. However
$f^j(x)=f^{j-s}(f^s(x))=f^{j-s}(0)$, and similarly for $f^j(y)$.  It
follows that if $f^j(X)+f^j(Y)$ fails to be absolutely irreducible,
then $f^{j-s}(0)+f^{j-t}(0)=0$ for some pair of non-negative integers
$s,t\le j-1$. If $s=t$ then since $\F_q$ has odd characteristic we
have $f^{j-s}(0)=0=f^0(0)$ with $1\le j-s\le
j$.  Otherwise $f^{j-s+1}(0)=f^{j-t+1}(0)$ with distinct positive
integers $j-s+1,j-t+1\le j+1$.  Since Theorem 1 assumes that the
values $f^0(0),f^1(0),\ldots,f^r(0)$ are distinct we therefore
conclude that the polynomial $f^j(X)+f^j(Y)$ is irreducible over the
algebraic completion $\overline{\F_q}$, for every $j<r$.

We are now ready to estimate $N(r;2)$.  In view of (\ref{N2}) and
(\ref{fact}) we have
\[N(r;2)\le q+\sum_{j=0}^{r-1}\#\{(x,y)\in\F_q^2:f^j(x)+f^j(y)\},\]
there being $q$ solutions to $x-y=0$. To get a corresponding lower
bound we may use the inclusion-exclusion principle to show that
\[N(r;2)\ge q+\sum_{j=0}^{r-1}\#\{(x,y)\in\F_q^2:f^j(x)+f^j(y)\}-
\sum_{0\le j\le r-1}A_j-\sum_{0\le i<j\le r-1}B_{ij},\]
where $A_j$ is the number of common solutions to
\[X-Y=0\;\;\;\mbox{and}\;\;\;f^j(X)+f^j(Y)=0, \]
and $B_{ij}$ is the number of common solutions to
\[f^i(X)+f^i(Y)=0\;\;\;\mbox{and}\;\;\;f^j(X)+f^j(Y)=0.\]

However if $f^j(x)+f^j(y)=0$ with $x=y$ then $f^j(x)=0$, which has at
most $2^j$ solutions. Thus $A_j\le 2^j$.  Similarly,
if $(x,y)$ were to lie on two distinct curves
$f^j(X)+f^j(Y)=0$ and $f^i(X)+f^i(Y)=0$ with $0\le i<j\le r-1$, then
\[f^j(y)=f^{j-i}(f^i(y))=f^{j-i}(-f^i(x))=f^j(x),\]
since $f^{j-i}$ is an even polynomial.  We would then have $2f^j(x)=0$
so that $x$, and similarly $y$, would be a root of $f^j$.  There are
therefore at most $2^j$ choices for $x$, and since $y$ then
satisfies $f^i(y)=-f^i(x)$ there are at most $2^i$ choices of $y$
for each possible $x$.  Thus $B_{ij}\le 2^{j+i}$.  It follows that
\[\sum_{0\le j\le r-1}A_j\le 2^r\;\;\;\mbox{and}\;\;\;
\sum_{0\le i<j\le r-1}B_{ij}\le 2^{2r}.\]
We therefore conclude that
\[N(r;2)= q+\sum_{j=0}^{r-1}\#\{(x,y)\in\F_q^2:f^j(x)+f^j(y)\}+O(4^r).\]

It remains to count points on the curves $f^j(X)+f^j(Y)=0$.  We have
already shown that these are absolutely irreducible, and indeed
nonsingular, under the assumptions of Theorem 1. If we write $N_r$ for
the number of projective points on the curve, and $D=2^j$ for its
degree, then Weil's ``Riemann Hypothesis'' tells us that
\[|N_r-(q+1)|\le (D-1)(D-2)\sqrt{q}.\]
There are at most $D$ points at infinity, so that
\[\left|\#\{(x,y)\in\F_q^2:f^j(x)+f^j(y)\}-q\right|\le D^2\sqrt{q}.\]
Finally, summing for $0\le j\le r-1$ we find that
\[\sum_{j=0}^{r-1}\#\{(x,y)\in\F_q^2:f^j(x)+f^j(y)\}=rq+O(4^r\sqrt{q}).\]

We may therefore summarize the conclusions of this section as follows.
\begin{lemma}
Under the assumptions of Theorem 1 we have
\[N(r;2):=\sum_{m\in \F_q}\rho_r(m)^2=(r+1)q+O(4^r\sqrt{q}).\]
\end{lemma}

\section{Higher Moments --- Irreducible Curves}

We now develop the ideas of the previous section to estimate $N(r;k)$
for $k\ge 3$. Here $N(r;k)$ is the number of solutions of
\beql{eqns}
f^r(x_1)=\ldots=f^r(x_k)
\eeq
in $\F_q$.  These equations define a curve, but, as in the previous
section, it is far from being an irreducible curve. Our task in this
section is to identify the absolutely irreducible components, and to
show that they are all defined over $F_q$.  

In view of (\ref{fact}), for any solution of (\ref{eqns}) and any pair
of distinct indices $1\le i,j\le k$, there is a
corresponding
\[d=d(i,j)=d(j,i)\in\{-1,0,1,\ldots,r-1\}\]
such that 
$\phi(x_i,x_j;d)=0$, where
\[\phi(X,Y;d)=\left\{\begin{array}{cc} f^d(X)+f^d(Y),& d\ge 0,\\
X-Y, & d=-1.\end{array}\right.\]
If there is more than one choice for $d(i,j)$ we choose the smallest.

We now make the following definition.
\begin{definition}
A ``$(D,k)$-graph'' is a weighted graph on $k$ vertices, for
which any edge $ij$ has integral weight in the range $[-1,D]$. If
some edge has weight equal to
$D$ we say that we have a ``strict $(D,k)$-graph''. If there
is an edge between every pair of vertices we say we have a ``complete
$(D,k)$-graph''.
\end{definition}
  
Thus each solution of (\ref{eqns}) produces a complete
$(D,k)$-weighted graph. We now introduce the following further
definition.

\begin{definition}\label{prop}
Let $G$ be a complete $(D,k)$-graph.  Then we say $G$ is 
 ``proper" if, whenever
$a,b,c$ are distinct vertices, with $d(a,b)\le d(a,c)\le d(b,c)$,
then either $d(a,b)=d(a,c)=d(b,c)=-1$ or
$d(a,b)<d(a,c)=d(b,c)$.  
\end{definition}

We then have the following lemma.
\begin{lemma}
  The graph associated to a solution of (\ref{eqns}) is proper.
\end{lemma}

To prove the claim, observe firstly that if 
$d(a,b)=d(a,c)=-1$, then $x_a=x_b$ and $x_a=x_c$, whence $x_a=x_c$,
so that $d(b,c)=-1$.  Next we show that one cannot have
$d(a,b)=d(a,c)\ge 0$.  Writing $d=d(a,b)=d(a,c)$ this would imply that
$f^d(x_a)=-f^d(x_b)$ and $f^d(x_a)=-f^d(x_c)$, whence
\[f^d(x_b)-f^d(x_c)=0.\]
The factorization (\ref{fact}) would then
show that $\phi(x_b,x_c;e)=0$ for some $e<d\le d(b,c)$. This however
is impossible, since $d(b,c)$ was chosen minimally.

To complete the proof of the claim we show that if
\[d(a,b)<d(a,c)\le d(b,c)\]
then $d(a,c)=d(b,c)$. In view of
(\ref{fact}) the relation
$\phi(x_a,x_b;d(a,b))=0$ would imply
$f^{d(a,c)}(X_a)-f^{d(a,c)}(X_b)=0$. Since
\[\phi(X_a,X_c;d(a,c))=f^{d(a,c)}(X_a)+f^{d(a,c)}(X_c)=0\]
this would show that
$f^{d(a,c)}(X_b)+f^{d(a,c)}(X_c)=0$ and the minimal choice of $d(b,c)$ 
then produces $d(b,c)\le d(a,c)$, giving the required conclusion.
This now establishes the lemma in full.

Thus each solution of (\ref{eqns}) is associated to a unique
proper weighted graph, such that
\beql{eqns1}
\phi(x_i,x_j;d(i,j))=0\;\;\;(1\le i\not= j\le k).
\eeq
However there is considerable redundancy in the equations
(\ref{eqns1}). To investigate this we begin with the following result.

\begin{lemma}\label{split}
  Let $G$ be a proper strict $(D,k)$-graph, with $D\ge 0$. Then there is a
unique partition $\{1,\ldots,k\}=A\cup B$ into non-empty sets $A$ and
$B$ such that $d(a,b)=D$ when $a\in A$ and $b\in B$, while $d(i,j)<D$
whenever $i,j\in A$ or $i,j\in B$.
\end{lemma}

Firstly it is easy to see that such a partition must be unique.  For
if $A'\cup B'$ were a different partition then, after relabeling if
necessary, we could find indices $i,j\in A$ with $i\in A'$ and $j\in
B'$.  We would then have both $d(i,j)<D$ (because $i,j\in A$) and
$d(i,j)=D$ (because $i\in A'$ and $j\in B'$). This contradiction shows
that such partitions are unique.

In order to show the existence of a suitable partition we fix a pair
$i_0,j_0$ with $d(i_0,j_0)=D$, and let
\[A=\cup\{i:d(i,j_0)=D\},\;\;\; B=\{j:d(j,i_0)=D\}.\]
Then $i_0\in A$ and $j_0\in B$, so that neither set is empty. If $a$, say,
were in $A\cap B$, then
\[d(a,j_0)=d(a,i_0)=d(i_0,j_0)=D\ge 0,\]
contradicting Definition \ref{prop}.
For any $a\in\{1,\ldots,k\}$ Definition \ref{prop}
shows that we must have either $d(a,i_0)=D$ or $d(a,j_0)=D$, so
that $A\cup B$ is a partition of $\{1,\ldots,k\}$. 

If $a_1,a_2\in A$ had $d(a_1,a_2)=D$ then the triple $a_1,a_2,j_0$
would contradict Definition \ref{prop}.  Thus
$d(a_1,a_2)<D$, and similarly $d(b_1,b_2)<D$ when $b_1,b_2\in B$.
Finally, if $a_1\in A$ then $d(a_1,j_0)=D$ if $a_1=i_0$.  Otherwise
Definition \ref{prop} applied to the triple $a_1,i_0,j_0$
shows that $d(a_1,j_0)=D$, since $d(a_1,i_0)<D$.  Thus $d(a,j_0)=D$
for all $a\in A$.  Now, if $b\in B$ with $b\not=j_0$, Definition
\ref{prop} applied to the triple $a,b,j_0$ shows that
$d(a,b)=D$, since $d(b,j_0)<D$.  Hence $d(a,b)=D$ whenever $a\in A$
and $b\in B$. This completes the proof of the lemma.

We now show how a complete $(D,k)$-graph can be generated by a smaller graph.

\begin{definition}\label{completion}
Let $G$ be a complete $(D,k)$-graph, and suppose $G_0$ is
a subgraph of $G$ with the same set of vertices but fewer edges.
We then say that $G_0$ ``generates'' $G$
if $G=G_n$ for some $n$, where $G_{h+1}$ is obtained from $G_r$
by the following procedure:

Take three distinct vertices $a,b,c$ for which the edges $ab$
and $bc$ belong to $G_r$ but $ac$ does not, and for which either
$d(a,b)=d(b,c)=-1$ or
$d(a,b)<d(b,c)$.  Then $G_{h+1}$ is obtained from $G_h$ by adding the
edge $ac$ with weight $d(a,c)=d(b,c)$.
\end{definition}

For our purposes it is not necessary to know whether, using a
different sequence of edge additions, $G_0$ might generate two
different complete $(D,k)$-graphs. All we need to know is whether
there exist some sequence of edge additions resulting in $G$.

To motivate the definition we consider the ideal
$I_h\subseteq\F_q[X_1,\ldots,X_k]$ generated by those polynomials
$\phi(X_i,X_j;d(i,j))$ for which the edge $ij$ is in $G_h$.  Then
trivially we have $I_h\subseteq I_{h+1}$, since $I_{h+1}$ is formed
from $I_h$ by the addition of one further generator
$\phi(X_a,X_c,d(a,c))$.  However, if $d(a,b)=d(b,c)=-1$ in the
procedure in Definition \ref{completion} we have
\[\phi(X_a,X_b,d(a,b))=X_a-X_b\]
and
\[\phi(X_b,X_c,d(b,c))=X_b-X_c.\]
Hence if $d(a,c)=-1$ then
\[\phi(X_a,X_c,d(a,c))=X_a-X_c=\phi(X_a,X_b,d(a,b))+\phi(X_b,X_c,d(b,c)),\]
so that $I_{h+1}=I_h$.  Alternatively, if
$d(a,c)=d(b,c)>d(a,b)$ in the
procedure in Definition \ref{completion}, we have
\begin{eqnarray*}
\phi(X_a,X_c,d(a,c))&=&\phi(X_a,X_c,d(b,c))\\
&=&f^{d(b,c)}(X_a)+f^{d(b,c)}(X_c)\\
&=&\left(f^{d(b,c)}(X_a)-f^{d(b,c)}(X_b)\right)+\phi(X_b,X_c:d(b,c)).
\end{eqnarray*}
Here we have $\phi(X_a,X_b;d(a,b))\mid
f^{d(b,c)}(X_a)-f^{d(b,c)}(X_b)$ by (\ref{fact}), since $d(a,b)<d(b,c)$.
Hence $\phi(X_a,X_c,d(a,c))$ is in the ideal generated by
$\phi(X_a,X_b,d(a,b))$ and $\phi(X_b,X_c,d(b,c))$.  We therefore see
again that $I_{h+1}=I_h$. It follows that
if $G$ is the proper complete $(D,k)$-graph associated to a system of
equation (\ref{eqns1}), and $G$ is generated by $G_0$, then the system
(\ref{eqns1}) has the same solutions as the smaller system
\beql{eqns2}
\phi(x_i,x_j;d(i,j))=0\;\;\;(ij\mbox{ is an edge of }G_0).
\eeq

We now introduce the small graphs we shall use.
\begin{definition}\label{chain}
A $(D,k)$-graph is said to be a ``chain'' if there is a permutation
$\sigma\in S_k$ such that the edges are precisely the $k-1$ pairs
\[(\sigma(1),\sigma(2))\,,\,(\sigma(2),\sigma(3))\,,\dots,\,
(\sigma(k-1),\sigma(k)),\]
and, for any $s<t\le k-1$, the maximum of
\[d(i_s,i_{s+1})\,,\,d(i_{s+1},i_{s+2})\,,\ldots,\,d(i_t,i_{t+1})\]
is either $-1$ or is attained at only one point.
\end{definition}

We then have the following result.
\begin{lemma}
For any complete $(D,k)$-graph $G$ there is a chain $(D,k)$-graph
$G_0$ which generates $G$.
\end{lemma}

We prove this by induction on $D$.  If $D=-1$ we may take $G_0$ to
consists of the edges $(1,2),\ldots,(k-1,k)$ with weights $-1$, which
clearly generates $G$.  Now assume the result is true for
complete $(d,k)$ graphs with $d\le D-1$.
If $G$ is not a strict $(D,k)$-graph the
conclusion is immediate from the induction hypothesis.

Thus we assume that $G$ is a strict complete $(D,k)$ graph with $D\ge
0$, so that Lemma \ref{split} applies. Let $G_A$ be the restriction of
$G$ to the vertices in $A$, so that $G_A$ is a complete
$(D-1,m)$-graph, where $m=\# A$.  The induction hypothesis then shows
that there is a chain graph $G_1$ say, which generates $G_A$, in which
one re-orders the vertices in $A$ as $i_1,\ldots,i_m$ so as to
satisfy the chain property in Definition \ref{chain}.
Similarly if $G_B$ is the restriction of
$G$ to the vertices in $B$, we can obtain
a subgraph $G_2$ of $G_B$ which is a chain, and which generates $G_B$. 
If $n=\# B$ there will again be an appropriate ordering
$j_1,\ldots,j_n$ of the indices in $B$.

We then take $G_0$ to be the graph with vertices $1,\ldots,k$ whose
edges are the edges of $G_A$, the edges of $G_B$, and the additional
edge $i_m,j_1$.  Moreover we permute the vertices into the order
$i_1,\ldots,i_m,j_1,\ldots,j_n$. We claim firstly that this ordering
makes $G_0$ a chain, and secondly that $G_0$ generates $G$.

To verify that $G_0$ is a chain we consider a sequence of consecutive
pairs of the vertices from the sequence $i_1,\ldots,i_m,j_1,\ldots,j_n$.
If the sequence is entirely contained in the first $m$ terms the
required chain property follows from that for $G_1$, and similarly if
all the elements are taken from the last $n$ terms.  However if one of
the pairs is the edge $i_m,j_1$ it suffices to note that this edge has
weight $D$ while all other edges have weight at most $D-1$.

To check that $G_0$ generates $G$ we note that $G_1$ generates $G_A$
and $G_2$ generates $G_B$.  Thus $G_0$ certainly generates the graph
$G^*$ containing the edges of $G_A$, the edges of $G_B$ and the edge
$i_m,j_1$. Hence it suffices to show that $G_*$ generates $G$.
Let $ij$ be an edge of $G$ which is not
already an edge in $G^*$. Then, according to Lemma
\ref{split} we may assume that $i\in A$ and $j\in B$, and that $d(i,j)=D$.
Applying the procedure in Definition \ref{completion} to the triple
$i,i_m,j_1$ we see that the edge $i,i_m$ is in $G^*$, since $i,i_m\in
A$, and the edge $i_m,j_0$ is also in $G^*$, by definition.  Moreover
$d(i,i_m)<D=d(i_m,j_1)$.  Thus the edge $i,j_1$ can be generated from
$G^*$, with weight $d(i,j_1)=d(i_m,j_1)=D$.  We may then apply
the procedure in Definition \ref{completion} to the triple
$j,j_1,i$. This time the edge $j,j_1$ is in $G^*$, since $j,j_1\in
B$, and the edge $j_1,i$ can be generated from $G^*$, as we have just
shown. Moreover we have $d(j,j_1)<D=d(j_1,i)$, so that the edge $ji$
can also be generated from $G^*$, and is given weight $D$, as
required. This completes the proof of the lemma.

As an immediate consequence of Lemma \ref{chain} we have the following.
\begin{lemma}
After a suitable relabelling of the variables, any solution to the equations
(\ref{eqns})
satisfies some system of equations of the type
\[\phi(X_i,X_{i+1};d_i)=0\;\;\; (1\le i\le k-1)\]
with $-1\le d_i\le k-1$.  Moreover, if
$1\le i<j\le k-1$, then the maximum of $d_i,\ldots, d_j$
is either $-1$ or occurs at only one point.
\end{lemma}

We call a system of equations of the above type a ``chain
system''. The system defines a variety in $\mathbb{A}^k$.  We set
$\Phi(X,Y,Z;-1)=X-Y$ and
\[\Phi(X,Y,Z;d)=Z^{2^d}\phi(X/X,Y/Z;d),\;\;(d\ge 0)\]
so that the corresponding projective variety is given by
\[\mathcal{C}:\,\Phi(X_i,X_{i+1},X_0;d_i)=0\;\;\; (1\le i\le k-1).\]

The importance of the chain property is demonstrated by the following
result.

\begin{lemma}\label{L4}
 Suppose that $f^i(0)\not=f^j(0)$ for $0\le i<j\le r$. Then,
  for a chain system, the variety $\mathcal{C}$ is a nonsingular
  complete intersection.  Hence $\mathcal{C}$
  is an absolutely irreducible curve
over $\F_q$, with degree at most $2^{(k-1)(r-1)}$.
\end{lemma}

To prove that $\mathcal{C}$ is a nonsingular complete intersection
we need to show that the vectors $\nabla\Phi(x_i,x_{i+1},x_0;d_i)$ are
linearly independent at any point of $\mathcal{C}$.  Suppose to the
contrary that
\[\sum_{i=1}^{k-1}c_i\nabla\Phi(x_i,x_{i+1},x_0;d_i)=\mathbf{0}.\]
If the $c_i$ are not all zero we take $s$ to be the smallest index
with $c_s\not=0$, and $t$ to be the largest index with $c_t\not=0$, so
that
\[\sum_{i=s}^tc_i\nabla\Phi(x_i,x_{i+1},x_0;d_i)=\mathbf{0}.\]
The entries of this vector are labelled by the variables
$X_0,\ldots,X_k$, and one sees that the entry corresponding to $X_s$
is just $c_s(\partial/\partial x_s)\Phi(x_s,x_{s+1},x_0;d_s)$.  We
therefore conclude that $(\partial/\partial
x_s)\Phi(x_s,x_{s+1},x_0;d_s)=0$, and similarly that $(\partial/\partial
x_{t+1})\Phi(x_t,x_{t+1},x_0;d_t)=0$. In particular we must have
$d_s,d_t\ge 1$.  However
\[\frac{\partial}{\partial x}\Phi(x,y,Z;d)=
(2a)^d\prod_{i=0}^{d-1}F^i(x,Z)\;\;\;(d\ge 0)\]
in the notation (\ref{nn}). We therefore see that $F^i(x_s,x_0)=0$ for
some index $i$ in the range $0\le i\le d_s-1$, and similarly
$F^j(x_{t+1},x_0)=0$ for some $j$ with $0\le j\le d_t-1$.

We next show that $x_0$ cannot vanish.  If, on the contrary, we had
$x_0=0$ then the relation $F^i(x_s,x_0)=0$ would yield $x_s=0$.  In
general, if $x_i=x_0=0$, then the relation $\Phi(x_i,,x_{i+1},x_0;d_i)=0$
implies $x_{i+1}=0$, while $\Phi(x_{i-1},,x_i,x_0;d_{i-1})=0$ implies
$x_{i-1}=0$. Thus, using both forwards and backwards induction we would
have $x_i=0$ for all $i$, which is impossible.

We may therefore assume that $x_0=1$, taking us back to the affine
situation. Thus we have $f^i(x_s)=0$ and $f^j(x_{t+1})=0$ with $0\le
i<d_s$ and $0\le j<d_t$. Since $d_s\ge 1$ the chain property shows
that the maximum of $d_s,d_{s+1},\dots,d_t$ occurs at only one point,
$d_u=D$, say. Since $i<d_s\le D$ we have
$f^D(x_s)=f^{D-i}(f^i(x_s))=f^{D-i}(0)$. Similarly we have
$f^D(x_{t+1})=f^{D-j}(0)$. If $s\le h<u$ then $d_h<D$, whence
$\phi(X,Y;d_h)\mid f^D(X)-f^D(Y)$.  Thus $f^D(x_h)=f^D(x_{h+1})$ for
$s\le h<u$.  It follows that $f^D(x_u)=f^D(x_s)=f^{D-i}(0)$.
Similarly, when $u<h\le t$ we have $f^D(x_h)=f^D(x_{h+1})$, whence
$f^D(x_{u+1})=f^D(x_{t+1})=f^{D-j}(0)$. However
$\phi(x_u,x_{u+1};D)=0$ with $D\ge 1$, whence
$f^D(x_u)+f^D(x_{u+1})=0$.  As in the previous section we
therefore conclude that
$f^{D-i}(0)+f^{D-j}(0)=0$ for some pair of non-negative integers
$i,j<D$. This leads either to $f^{D-i}(0)=0$ (if $i=j$) or
$f^D(0)=f^{D+i-j}(0)$ (if $i<j$, say). In either case we contradict the
assumption of Theorem 1, since $D\le r$.  This completes the proof
that $\mathcal{C}$ is a nonsingular complete intersection.

The remainder of the lemma is now straightforward.  In general a
nonsingular complete intersection is necessarily absolutely
irreducible, with degree equal to the product of the degrees of the
defining forms, see Browning and Heath-Brown \cite[Lemma 3.2]{many} for
details.  In our case $\Phi(X_i,X_{i+1},X_0;d_i)$ has degree at most
$2^{r-1}$, since $d_i\le r-1$, and the result follows.

\section{Higher Moments --- Counting Points, And Counting Curves}

In this section we will firstly estimate the number of points on each curve
$\mathcal{C}$, and then compute the number of such curves that the
variety given by (\ref{eqns}) produces.  Putting these results
together will give us an asymptotic formula for $N(r;k)$.

Since $\mathcal{C}$ is an absolutely irreducible curve defined over
$\F_q$, Weil's ``Riemann Hypothesis'' yields
\[\left|\#\mathcal{C}(\F_q)-(q+1)\right|\le 2g\sqrt{q},\]
where $g$ is the genus of $\mathcal{C}$.  In general, if $\mathcal{C}$
is an irreducible non-degenerate curve of degree $d$ in
$\mathbb{P}^k$ (with $k\ge 2$), then
according to the Castelnuovo genus bound \cite{Cast}, one has
\[g\le (k-1)m(m-1)/2+m\varepsilon,\]
where $d-1=m(k-1)+\varepsilon$ with $0\le\varepsilon<k-1$. This
implies in particular that $g\le (d-1)(d-2)/2$ irrespective of the
degree of the ambient space in which $\mathcal{C}$ lies.  We
therefore deduce that
\beql{wb}
\left|\#\mathcal{C}(\F_q)-(q+1)\right|\le 4^{kr}\sqrt{q},
\eeq
since $\mathcal{C}$ has degree at most $2^{kr}$.

By inclusion-exclusion we see that
\[\sum_{\mathcal{C}}\#\mathcal{C}(\F_q)-
\frac{1}{2}\sum_{\mathcal{C}_1\not=\mathcal{C}_2}
\#\left(\mathcal{C}_1\cap\mathcal{C}_2\right)(\F_q)\leq N(r;k)\le
\sum_{\mathcal{C}}\#\mathcal{C}(\F_q).\]
For distinct curves of degree at most $2^{kr}$ we have
\[\#\left(\mathcal{C}_1\cap\mathcal{C}_2\right)(\F_q)\leq 4^{kr},\]
by B\'{e}zout's Theorem.  Hence if there are $\mathcal{N}(r;k)$
different curves $\mathcal{C}$ we see that
\[\left|N(r;k)-\sum_{\mathcal{C}}\#\mathcal{C}(\F_q)\right|
\leq \mathcal{N}(r;k)^2 2^{k2^r}.\]
We can get a crude bound for $\mathcal{N}(r;k)$ by observing that
there are $k!$ possible permutations describing a chain system, and
for each of the $k-1$ edges one has $-1\le
d(\sigma(i),\sigma(i+1))\le r-1$.  Thus if $r\ge 1$ we have
\beql{Nb}
\mathcal{N}(r;k)\le k!(r+1)^{k-1}\le k!(2r)^{k-1}\le (rk)^k\le
(2^rk)^k. 
\eeq
Applying (\ref{wb}) we then deduce the following result.

\begin{lemma}\label{NL}
If there are $\mathcal{N}(r;k)$ different curves $\mathcal{C}$ then
\[N(r;k)=\mathcal{N}(r;k)(q+1)+O(2^{4kr}k^{2k}\sqrt{q}).\]
\end{lemma}

Our task now is to investigate the number $\mathcal{N}(r;k)$. We have
seen that each curve $\mathcal{C}$ arises from a proper
$(r-1,k)$-graph.  We proceed to show that different graphs $G,G'$
cannot produce the same curve $\mathcal{C}$.  The graphs $G$ and $G'$
must differ on at least one
edge, so that one would have both $\phi(X_i,X_j;d)=0$ and
and $\phi(X_i,X_j;d')=0$ on $\mathcal{C}$. If $d>d'$ say, then
$f^d(X_i)+f^d(X_j)=0$ and $f^d(X_i)-f^d(X_j)=0$, whence
$f^d(X_i)=0$ for all points on the curve. It then follows that
$f^r(X_i)=f^{r-d}(0)$. However $\mathcal{C}$ is an irreducible
component of the curve (\ref{eqns}), whence
$f^r(X_h)=f^r(X_i)=f^{r-d}(0)$ for every index $h$. This gives us a
contradiction since it would produce imply that $\mathcal{C}$ has
dimension zero.

We therefore need to count proper $(r-1,k)$ graphs. For a proper
strict $(D,k)$-graph, Lemma \ref{split} produces a unique partition
$A\cup B$, for which the corresponding graphs $G_A$ and $G_B$ will be
proper $(D-1,k)$-graphs.  There are
$\mathcal{N}(r;k)-\mathcal{N}(r-1;k)$ proper strict $(r-1,k)$-graphs.
Moreover the number of partitions $\{1,\ldots,k\}=A\cup B$ with
$a=\# A<\# B=b$ is
\[\left(\begin{array}{cc} k\\ a\end{array}\right),\]
  while, for even $k$, the number with $a=b=k/2$ is
\[\frac{1}{2}\left(\begin{array}{cc} k\\ k/2\end{array}\right).\]
We then see that
\[\mathcal{N}(r;k)-\mathcal{N}(r-1;k)=
\frac{1}{2}\sum_{a=1}^{k-1}\left(\begin{array}{cc}
  k\\ a\end{array}\right) \mathcal{N}(r-1;a)\mathcal{N}(r-1;k-a)\]
  for $r\ge 1$ and $k\ge 2$. Indeed, since $\mathcal{N}(r;1)=1$ for every
  $r$ we see that this holds for $k=1$ too.
If we now define $\mathcal{N}(r;0)=1$ for all $r\ge 0$ the above
formula simplifies to give
\[\mathcal{N}(r;k)=
\frac{1}{2}\sum_{a=0}^k\left(\begin{array}{cc}
  k\\ a\end{array}\right) \mathcal{N}(r-1;a)\mathcal{N}(r-1;k-a)
  \;\;\; (r,k\ge 1).\]

We therefore define power series
  \[E(X;r):=\sum_{k=0}^{\infty}  \frac{\mathcal{N}(r;k)}{k!}X^k,\]
  for each $r\ge 1$.  Since (\ref{Nb}) yields $\mathcal{N}(r;k)/k!\le
(r+1)^k$ we see that this converges absolutely for $|X|<(r+1)^{-1}$.
 Now, after checking that we have the correct
  coefficient for $X^0$, we arrive at
  \[E(X;r)^2=\frac{1+E(X;r-1)^2}{2}\;\;\;(r\ge 1).\]
  Since $\mathcal{N}(0;k)=1$ for all $k$ we have $E(X;0)=\exp(X)$, so that
  the coefficients $\mathcal{N}(r;k)$ can easily be calculated in
  general. Moreover it is clear by induction that
  \[E(X;r)=\sum_{m=0}^{2^r}\nu(r;m)e^{mX}\]
  with non-negative real coefficients $\nu(r;m)$ summing to 1.
  We then see that
  \[E(X;r)=\sum_{m=0}^{2^r}\nu(r;m)\sum_{k=0}^{\infty}\frac{(mX)^k}{k!}.\]
 We clearly have have absolute
convergence for small $X$, and we may rearrange to get
\[E(X;r)=
\sum_{k=0}^{\infty}\left(\sum_{m=0}^{2^r}\nu(r;m)m^k\right)\frac{X^k}{k!}.\]
 We therefore deduce that
 \[\mathcal{N}(r;k)=\sum_{m=0}^{2^r}\nu(r;m)m^k.\]
We also see that the coefficient $\nu(r;0)$ satisfies the recurrence
\[\nu(r;0)=\frac{1+\nu(r-1;0)^2}{2}\]
for $r\ge 1$, with 
$\nu(0;0)=0$. We can then check that $\mu_r=1-\nu(r;0)$ has the
initial value $\mu_0=1$ and satisfies the recurrence
$\mu_r=\mu_{r-1}-\mu_{r-1}^2/2$ described in Theorem 1.

Recall that our goal is to estimate
  \[\# f^r(\F_q)=q-\#\{m\in\F_q:\rho_r(m)=0\}.\]
  Since the equation $f^r(x)=m$ has at most $2^r$ solutions we will always
  have $0\le \rho_r(m)\le 2^r$, whence
  \[\frac{1}{2^r!}\prod_{j=1}^{2^r}\left(j-\rho_r(m)\right)=
  \left\{\begin{array}{cc}
    1, & \rho_r(m)=0,\\  0, & \rho_r(m)>0.\end{array}\right.\]
Setting
\beql{iden}
\frac{1}{2^r!}\prod_{j=1}^{2^r}\left(j-T\right)=\sum_{k=0}^{2^r}C_{r,k}T^k
\eeq
we then have
\beql{st}
\#\{m\in\F_q:\rho_r(m)=0\}=\sum_{k=0}^{2^r}C_{r,k}N(r;k).
\eeq
Our plan is to substitute the approximate value for $N(k;r)$ given by
Lemma \ref{NL}.

We first investigate the contribution from the main term
$\mathcal{N}(k,r)(q+1)$.  This produces
\begin{eqnarray*}
(q+1)\sum_{k=0}^{2^r}C_{r,k}\,\mathcal{N}(r;k)&=&
  (q+1)\sum_{k=0}^{2^r}C_{r,k}\sum_{m=0}^{2^r}\nu(r;m)m^k\\
  &=& (q+1)\sum_{m=0}^{2^r}\nu(r;m)\sum_{k=0}^{2^r}C_{r,k}m^k.
  \end{eqnarray*}
However the identity (\ref{iden}) shows that the inner sum vanishes for
$1\le m\le 2^r$, and takes the value 1 for $m=0$.  Thus the main term
for (\ref{st}) is just $\nu(r;0)(q+1)=(1-\mu_r)(q+1)$, producing
the leading term $\mu_rq$ in (\ref{stst}).

For the proof of Theorem 1 it remains to handle the contribution to
(\ref{st}) arising from the error term in Lemma \ref{NL}, which will
be
\[\ll\sqrt{q}\sum_{k=0}^{2^r}|C_{r,k}|4^{kr}k^{2k}
\le \sqrt{q}\sum_{k=0}^{2^r}|C_{r,k}|16^{kr}=
\sqrt{q}G_r(16^r),\] 
with
\[G_r(T):=\sum_{k=0}^{2^r}|C_{r,k}|T^k.\]
However it is clear from (\ref{iden}) that
\[G_r(T)\le \frac{1}{2^r!}\prod_{j=1}^{2^r}\left(j+T\right)
\le \max(2^r,T)^{2^r}\]
if $T\ge 0$, so that
\[G_r(16^r)\le\{16^r\}^{2^r}\ll 2^{4^r},\]
say.  This suffices for Theorem 1.

\bigskip

\bigskip

Mathematical Institute,

Radcliffe Observatory Quarter,

Woodstock Road,

Oxford

OX2 6GG

UK

\bigskip

{\tt rhb@maths.ox.ac.uk}

\end{document}